\documentclass[10pt]{article}

\title{\bf The periodic Cauchy problem of the modified Hunter-Saxton equation}

\author{Feride T\i\u{g}lay}

\usepackage{amsmath, amssymb, amsthm}

\def\eps{\varepsilon}
\def\ga{\gamma}
\def\gi{\gamma^{-1}}
\def\gd{\dot{\gamma}}

\def\z{\zeta}

\def\diff{\mathcal{D}^{s}}
\def\z{\zeta}

\def\zd{\dot{\zeta}}
\def\dt{\partial_{t}}
\def\dx{\partial_{x}}
\def\dxi{\partial_{x}^{-1}}
\def\dxg{(\partial_{x})_{\ga}}
\def\dxig{(\partial_{x}^{-1})_{\ga}}

\def\sob{H^{s}}
\def\sobb{H^{s-1}}
\def\sod{H^{s}}
\def\sodd{H^{s-1}}
\def\dif{\mathcal{D}^{s}}

\def\tor{\mathbb{T}}

\def\reel{\mathbb{R}}

\def\D{\mathcal{D}}

\newtheorem {defi}{Definition}
\newtheorem {lema} {Lemma}
\newtheorem {prop}{Proposition}
\def \pf {\textit{Proof. }}

\newtheorem {rema}{Remark}
\newtheorem {teor}{Theorem}

\begin{document}
\maketitle
\begin{abstract}
We prove that the periodic initial value problem for the modified Hunter-Saxton equation is locally well-posed for initial data in the space of continuously differentiable functions on the circle and in Sobolev spaces $\sob (\tor)$ when $s>3/2$. We also study the analytic regularity (both in space and time variables) of this problem and prove a Cauchy-Kowalevski type theorem. Our approach is to rewrite the equation and derive the estimates which permit application of o.d.e. techniques in Banach spaces. For the analytic regularity we use a contraction argument on an appropriate scale of Banach spaces to obtain analyticity in both time and space variables.
\end{abstract}

In this paper we study the periodic Cauchy problem of the modified Hunter-Saxton equation
\begin{equation}\begin{array}{l}
\dt u +u^{p}\dx u=\frac{1}{2} \dxi (\dx (u^{p})\dx u ), \\
u(x,0)=u_{0}(x), \ x\in\tor \simeq \reel / \mathbb{Z} , t\in\reel .
\end{array}
\tag{mHS}
\end{equation}
where $p$ is any positive integer.

J.K. Hunter and R. Saxton derived this family of equations in \cite{HS} and showed that smooth solutions on $\reel$ break down in finite time. 

For $p=1$ the equation 
\begin{equation}
\dt u + u \dx u = \dxi\frac{\textstyle 1}{\textstyle 2} (\dx u)^{2} 
\tag{HS}
\end{equation} 
is called Hunter-Saxton equation.

Written in the form (HS) the Hunter-Saxton equation can be viewed as a nonlocal perturbation of the Burgers equation like the Camassa-Holm equation. Both equations Camassa-Holm and Hunter-Saxton arise along with the Korteweg-De Vries equation in these nonlocal forms when derived on the Bott-Virasoro group \cite{KhMi}. Similar local well-posedness results hold for Camassa-Holm equation in Sobolev spaces (see for example the proof given by A. Himonas and G. Misio\l ek in \cite{HM2} or \cite{Mis1}) and in $C^{1}$ (the $C^{1}$ theory for Camassa-Holm equation is developed by G. Misio\l ek in \cite{Mis1}). Another context in which the three equations Korteweg-De Vries, Camassa-Holm and Hunter-Saxton are treated in a unified way is their scattering theory described by R. Beals, D.H. Sattinger and J. Szmigielski in \cite{BSS1}. All three equations are bihamiltonian and have the same symmetry group, namely Virasoro group as shown in \cite{KhMi} by B. Khesin and G. Misio\l ek. Moreover these equations correspond to the equations of the geodesic flow with respect to different right invariant Riemannian metrics on this group or on an associated homogeneous space \cite{KhMi}.

An observation of V. Arnold \cite{Arn}, that the initial value problem for the classical Euler equations of a perfect fluid
can be stated as a problem of finding geodesics on the group of volume preserving
diffeomorphisms, was used by D.G. Ebin and J. Marsden in \cite{EMa} who
developed the necessary functional analytic tools and established sharp local well-posedness results for the Euler
equations in a class of Sobolev spaces. 

There is no similar result for
$p\geq 2$ for the modified Hunter-Saxton equation. It turns out, however, that the method of rewriting the problem as an ordinary differential equation on
the group of diffeomorphisms (just as in the case of Euler equations) can be applied by introducing two dependent variables namely $\z$ and $\ga$ in (\ref{eq:*}).
We use this approach to develop an appropriate analytic framework for (mHS) and prove the following theorems.
\begin{teor}
Let $p\geq 1$ be any positive integer and $s>3/2$. Given the initial data $u_{0} \in \sod(\tor)$
the Cauchy problem for the equation (mHS) has a unique solution
\[ u \in C^{0}([0, T), \sod(\tor)) \cap  C^{1} ([0, T), \sodd (\tor))
\]
for some $T>0$ and the solution depends continuously on the initial data $u_{0}(x)$.
\label{th:2}
\end{teor}

Clearly the local well-posedness in $C^{1}(\tor)$ as stated in theorem \ref{th:C1} is a stronger result than the local well-posedness in $\sob(\tor)$ for $s>3/2$ as stated in theorem \ref{th:2}. 

\begin{teor} 
Let $p\geq 1$ be any positive integer.
For any $u_{0}\in C^{1}(\tor)$ the Cauchy problem for (mHS) has a unique solution
\[ u\in C^{0}([0,T), C^{1}(\tor))\cap C^{1} ([0, T), L^{\infty} (\tor))
\]
for some $T>0$ and the solution depends continuously on the initial data.
\label{th:C1}
\end{teor}

After the completion of this paper we learnt that for $p=1$ theorem \ref{th:2} has been proved independently by Z. Yin in \cite{Yi} using the semigroup theory of T. Kato (see the remark \ref{rem:2} below).  

Furthermore the two equations, the modified Hunter-Saxton (mHS) and the Camassa-Holm, do not have similar features only for low regularity data but also for high regularity data. 
The study of analytic
regularity of solutions of the Camassa-Holm equation by A. Himonas and G. Misio\l ek \cite{HM1} using an abstract Cauchy-Kowalevski
theorem led us to investigate the analytic regularity of the Cauchy problem for (mHS) and prove the following theorem.
\begin{teor}
If the initial data $u_{0}$ is an analytic function on $\tor$ then there exists
an $\eps>0$ and a unique solution $u$ of the Cauchy problem for the equation (mHS) that is analytic in both variables $x$
and $t$ on $\tor$ for all $t$ in $(-\eps, \eps)$.
\label{th:analy_reg}  
\end{teor}

We would like to note that the analyticity properties of the solutions to the Hunter-Saxton and Camassa-Holm equations are quite different from those of the Korteweg-De Vries equation whose solutions are analytic in the space variable for all time \cite{Tru} but are not analytic in the time variable \cite{KaM}.

Our approach in proving theorem \ref{th:analy_reg} is to rewrite the equation to use a contraction argument on an appropriate  scale of Banach spaces. This contraction argument was developed in the form of an abstract Cauchy-Kowalevski theorem by L.V. Ovsjannikov \cite{Ovs1,Ovs2}, F. Treves
\cite{Tre}, L. Nirenberg \cite{Nir}, T. Nishida \cite{Nis} and M.S. Baouendi and C.
Goulaouic
\cite{BG} among others and subsequently applied to the Euler and Navier-Stokes equations. 

\section{Well-posedness in $C^{1}$}

Many techniques in studying differential equations are based on Picard's contraction argument on Banach spaces. Here we develop the tools to use this argument for the Cauchy problem of (mHS). Let us first introduce the operator $\dxi$ to write the modified Hunter-Saxton equation in the form (mHS):
\[ \dxi f(x):=\int_{x_{0}}^{x} f(y)dy-\int_{\tor}\int_{x_{0}}^{x}f(y)dy dx 
\]
for some fixed $x_{0}$ in $\tor$. 

The Hunter-Saxton equation (HS) can be restated on the space $\D$ of $C^{1}(\tor)$ diffeomorphisms of $\tor$  as an ordinary 
differential equation. 
The idea comes from its derivation as a geodesic equation on the Bott-Virasoro group (see \cite{KhMi}). Similarly one can restate (mHS) as an ordinary differential equation on the product space $\D\times C^{1}(\tor)$ as follows.

Let $ \z $ be equal to $u \circ \ga $ where $ \ga$ is the flow generated by $u^{p}$. Then we 
obtain the set of equations
\begin{equation}
\begin{array}{l}
\gd = \z^{p}, \\
\zd = ( \dt (\z \circ \gi) + (\z \circ \gi)^{p} \dx (\z \circ \gi)) \circ \ga.
\end{array}
\label{eq:*}
\end{equation}

The two initial value problems, for (mHS) and for the equation in (\ref{eq:*}), are equivalent in the following sense.

\begin{prop}
A function $u \in C^{1}(\tor)$ is a solution to the Cauchy problem
\begin{equation}
\dt u +u^{p} \dx u=\frac{p}{2} \dxi (u^{p-1} (\dx u)^{2}) ;  \ u(x,0)=u_{0}(x)
\label{eq:mhs_ini}
\end{equation}
if and only if  $u$ can be written as $\z \circ \gi$ where $(\ga , \z)\in\D \times C^{1}(\tor)$ is a solution to
\begin{equation}
\begin{array}{l} 
\gd = \z^{p}, \\ \\
\zd =F(\ga, \z) := \frac{\textstyle p}{\textstyle 2} ( \dxi ( (\z \circ \gi)^{p-1} (\dx (\z \circ \gi))^{2} ) ) \circ \ga
\end{array}
\label{eq:**1}
\end{equation}
with initial data $\z (x,0)=u_{0}(x)$ and $\ga (x,0) = id_{x} $.
\label{pro:odeC1}
\end{prop}

Therefore it is sufficient to prove that the pair $(\z^{p},F(\ga,\z))$ defines a continuously differentiable vector field in a neighborhood of $(id,0)$ in $\D \times C^{1}(\tor)$. Then theorem \ref{th:C1} follows by the fundamental theorem on ordinary differential equations in Banach spaces with the observation that the smooth dependence on initial data in (\ref{eq:**1}) implies continuous dependence on initial data for (\ref{eq:mhs_ini}).
\\
\\
{\bf \textit{Proof of Theorem \ref{th:C1}.}} Let us use the convenient notation from \cite{Eb1} and denote by $P_{\ga}$ the operator given by conjugation
\[ P_{\ga}(g):=P(g\circ \gi)\circ\ga
\]
for any $\ga\in \D$ and pseudodifferential operator $P$. Using this notation we write the right hand side of the first equation in (\ref{eq:**1}) as
\[ F(\ga,\z)=\frac{1}{2}\dxig g(\ga,\z)
\]
where $g(\ga,\z)=\dxg \z^{p}\dxg \z$.

Next we compute the directional derivatives $\partial_{\ga}F_{(\ga,\z)}$ and $\partial_{\z}F_{(\ga,\z)}$ and prove that these are bounded linear maps.

Note that it is enough to determine $\partial_{\ga}\dxig$ and $\partial_{\ga}\dxg$ since we have the following identities
\begin{equation}
\partial_{\ga}F_{(\ga,\z)}=\frac{1}{2} \partial_{\ga}\dxig g(\ga,\z)+\frac{1}{2}\dxig \partial_{\ga}g_{(\ga,\z)},
\label{eq:fga}
\end{equation}
\begin{equation}
\partial_{\z}F_{(\ga,\z)}=\frac{1}{2}\dxig \partial_{\z}g_{(\ga,\z)}
\label{eq:fz}
\end{equation}
and
\begin{equation}
\partial_{\ga}g_{(\ga,\z)}=\partial_{\ga}\dxg (\z^{p})\dxg \z +\dxg (\z^{p})\partial_{\ga}\dxg \z.
\label{eq:gg}
\end{equation}

Let $s\rightarrow \ga_{s}$ be a smooth curve in $\D$ such that $\ga_{0}=id$ and $\partial_{s} \ga_{s}|_{s=0}=W$ for $W \in C^{1}(\tor)$. By the definition of the operator $\dxi$ we obtain
\[ \partial_{\ga}\dxig (g)(W)=-\dxig (W \dxg g)+gW
\]
\begin{equation} 
= \dxig (g \dxg W)+ \int_{\tor}gW \ dx.
\label{eq:dgdxi}
\end{equation}
Similarly we have
\begin{equation} 
\partial_{\ga}\dxg =[W\circ \gi \dx , \dx ]_{\ga}
\label{eq:dgdx}
\end{equation}
where the bracket $[\cdot , \cdot]$ denotes the commutator of the operators.

By (\ref{eq:fga})-(\ref{eq:dgdx}) the directional derivatives of $F(\ga,\z)$ are given by the following explicit formulas:
\begin{eqnarray}
\partial_{\ga}F_{(\ga,\z)}(W)=-\frac{1}{2}\dxig (\dxg \z^{p}\dxg \z \dxg W) \nonumber \\
+\frac{1}{2}\int_{\tor} W \dxg \z^{p}\dxg \z \ dx,
\label{eq:dgf}  
\end{eqnarray}
\begin{equation}
\partial_{\z}F_{(\ga,\z)}(W)=\frac{1}{2}\dxig \{ p\dxg (\z^{p-1}W)\dxg \z + \dxg \z^{p}\dxg W \}.
\label{eq:dzf}  
\end{equation}

The linearity of the maps $W\rightarrow \partial_{\ga}F_{(\ga,\z)}(W)$ and $W\rightarrow \partial_{\z}F_{(\ga,\z)}(W)$ is obvious. Thus we proceed to show that these maps are bounded. For the $C^{1}$ norm of $\partial_{\ga}F_{(\ga,\z)}$ 
we have
\begin{eqnarray}
\| \partial_{\ga}F_{(\ga,\z)}(W)\|_{C^{1}} &\leq & \| \dxig (\dxg (\z^{p})\dxg\z \dxg W)\|_{C^{1}} \nonumber \\
& &+\left| \int_{\tor} W \dxg (\z^{p})\dxg \z \ dx\right| \nonumber \\
&\leq & C_{\ga}\| \dxg (\z^{p})\dxg\z \dxg W\|_{\infty}  \nonumber \\
& &+ \| \dxig (\dxg (\z^{p})\dxg\z \dxg W)\|_{\infty} \label{eq:bfg}\\
& &+ C_{\ga} \| W\|_{\infty} \| \z\|_{C^{1}}^{p+1} \nonumber
\end{eqnarray}
where $C_{\ga}$ depends only on $\|\ga\|_{C^{1}}$ and $\| \gi\|_{C^{1}}$.
The first term in (\ref{eq:bfg}) is bounded by
\[ C_{\ga} \| \z\|_{C^{1}}^{p+1} \| W\|_{C^{1}}.
\] 
We also observe that the second term in (\ref{eq:bfg}) can be written as
\[ \left\| \int_{\ga(x_{0})}^{\ga(x)} \dx (\z^{p}\circ \gi) \dx (\z\circ\gi)\dx (W\circ \gi)dy \right\|_{\infty}
\]
\[ + \left| \int_{\tor}\int_{\ga(x_{0})}^{\ga(x)} \dx (\z^{p}\circ \gi) \dx (\z\circ\gi)\dx (W\circ \gi)dy \ dx \right|.
\]
This last sum is bounded by
\[ C_{\ga} \| W\|_{C^{1}} \|\z\|_{C^{1}}^{p+1}. \]
Combining these estimates we obtain
\[ \|\partial_{\ga}F_{(\ga,\z)}(W)\|_{C^{1}}\leq C_{\ga} \|W\|_{C^{1}} \|\z\|_{C^{1}}^{p+1}
\]
and therefore the map $W\rightarrow \partial_{\ga}F_{(\ga,\z)}(W)$ is bounded.

The map $W\rightarrow \partial_{\z}F_{(\ga,\z)}(W)$ is bounded as well by the estimate 
\begin{eqnarray}
\| \partial_{\z}F_{(\ga,\z)}(W)\|_{C^{1}} & \leq & C_{\ga} \left( \| \dxg (W\z^{p-1})\dxg \z\|_{\infty}+\|\dxg\z^{p}\dxg W \|_{\infty}\right) \nonumber \\
& & +\| \dxig \{p\dxg (W\z^{p-1})\dxg \z+ \dxg \z^{p}\dxg W \}\|_{\infty} \nonumber \\
& \leq & C_{\ga} \| W\|_{C^{1}} \| \z\|_{C^{1}}^{p}. \nonumber
\end{eqnarray}
which is derived similarly.

In order to complete the proof of theorem \ref{th:C1} we show that both directional derivatives $\partial_{\ga}F_{(\ga,\z)}$ and $\partial_{\z}F_{(\ga,\z)}$ are continuous maps. Note that to prove the continuity of $\partial_{\ga}F_{(\ga,\z)}$ in $\ga$ (uniformly in $\z$) from $\D$ into the space of bounded linear operators on $C^{1}(\tor)$ it suffices to prove the continuity of the following maps
\begin{eqnarray}
 & & \ga\rightarrow\partial_{\ga}\dxig\in L(C^{1},L(C^{0},C^{1})) \label{eq:ayva} \\
 & & \ga\rightarrow\dxig\in L(C^{0},C^{1}) \label{eq:elma} \\
 & & \ga\rightarrow\dxg\in L(C^{1},C^{0}) \label{eq:nar}
\end{eqnarray}
in a neighborhood of $(id,0)\in \D\times C^{1}(\tor)$.
We consider the map (\ref{eq:elma}) first. A change of variables in
\begin{eqnarray*}
\dxig f- \dxi f = \int_{\ga(x_{0})}^{\ga(x)}f\circ\gi(y)dy- \int_{\tor}\int_{\ga(x_{0})}^{\ga(x)}f\circ\gi(y)dy \ dx \\
- \int_{x_{0}}^{x}f(y)dy+ \int_{\tor}\int_{x_{0}}^{x}f(y)dy \ dx 
\end{eqnarray*}
leads to the identity
\begin{equation}
\dxig f- \dxi f = \dxi (f(\dx\ga-1)).
\label{eq:identity}
\end{equation}
The right hand side of this identity is bounded by
\[ C_{\ga} \| f \|_{\infty} \| \ga- id\|_{C^{1}}
\]
and therefore we have the estimate
\begin{equation}
\|\dxig f- \dxi f\|_{C^{1}} \leq C_{\ga} \| f\|_{\infty} \| \ga- id\|_{C^{1}} 
\label{eq:ikiyildiz}
\end{equation}
and the map in (\ref{eq:elma}) is continuous.

In order to show the continuity of the map in (\ref{eq:ayva}) we observe that the definite integral in $\partial_{\ga}\dxig (W)$ does not depend on $x$ and thus we have
\begin{eqnarray}
\lefteqn{\| \partial_{\ga}\dxig (W)-\partial_{\ga}(\dxi)_{id}(W)\|_{C^{1}}} \nonumber \\
& \leq & \| \dxig \{ \dxg \z^{p}\dxg \z\dxg W \} - \dxi (\dx \z^{p}\dx\z\dx W)\|_{C^{1}} \label{eq:4star1} \\
& & +\left| \int_{\tor}W\{ \dxg \z^{p}\dxg \z-\dx\z^{p}\dx\z\}dx \right| \label{eq:4star2}
\end{eqnarray}  
Adding and subtracting the appropriate terms we can bound (\ref{eq:4star1}) by
\begin{eqnarray}
& & C_{\ga} \| \dx\z^{p}\dx\z \dx W \{ (\dx\gi\circ\ga)^{3}-1 \}\|_{\infty} \nonumber \\
& & +\| \dxig \{ \dx\z^{p}\dx\z \dx W ( (\dx\gi\circ\ga)^{3}-1 )\}\|_{\infty} \label{eq:ucyildiz}\\
& & +\| \{ \dxig - (\dxi)_{id}\}(\dx\z^{p}\dx\z\dx W)\|_{C^{1}}. \nonumber
\end{eqnarray}
We use (\ref{eq:ikiyildiz}) to estimate the third summand above by
\[ \| \z\|_{C^{1}}^{p+1}\| W\|_{C^{1}} \| \ga -id\|_{C^{1}}.
\]
Clearly the first two terms  in (\ref{eq:ucyildiz}) are bounded by 
\[ C_{\ga} \| \z\|_{C^{1}}^{p+1}\| W\|_{C^{1}} \| \ga -id\|_{C^{1}}
\]
and therefore we have the following estimate for (\ref{eq:4star1})
\[ C_{\ga} \| \z\|_{C^{1}}^{p+1}\| W\|_{C^{1}} \| \ga -id\|_{C^{1}}.
\]

We also observe that (\ref{eq:4star2}) can be bounded easily as follows
\begin{eqnarray*}
\lefteqn{\left| \int_{\tor}W\{ \dxg \z^{p}\dxg \z-\dx\z^{p}\dx\z\}dx \right|} \\
& \leq & \| W\|_{\infty} \left| \int_{\tor} \dx \z^{p}\dx \z ((\dx\gi\circ \ga)^{2}-1)dx \right| \\
& \leq & C_{\ga} \| W\|_{\infty} \|\z\|_{C^{1}}^{p+1}\|\ga-id\|_{C^{1}}.
\end{eqnarray*}
Combining this last estimate with the estimate for (\ref{eq:ucyildiz}) we obtain
\[ \| \partial_{\ga}\dxig (W)-\partial_{\ga}(\dxi)_{id} (W)\|_{C^{1}} \leq C_{\ga}\|W\|_{C^{1}}\|\z\|_{C^{1}}^{p+1}\|\ga-id\|_{C^{1}}
\]
which implies continuity of the map (\ref{eq:ayva}).

For the continuity of the map (\ref{eq:nar}) it is sufficient to observe that
\begin{eqnarray*}
\| \dx(f\circ\gi)\circ\ga-\dx f\|_{\infty}&=& \| \dx f (\dx\gi\circ\ga -1)\|_{\infty} \\
&\leq & C_{\ga} \|f\|_{C^{1}} \| \ga -id\|_{C^{1}}.
\end{eqnarray*}
Hence the continuity of $\ga \rightarrow \partial F_{\ga,\z}$ follows.

Furthermore the continuity in $\z$ of $\partial F_{\ga,\z}$ is easier to prove. 
We rewrite the $C^{1}$ norm that we want to estimate as
\begin{eqnarray}
\lefteqn{\| \partial_{\ga}F_{\ga,\z_{1}}-\partial_{\ga}F_{\ga,\z_{2}}\|_{C^{1}}} \nonumber \\
& \leq & \| \dxig \left\{ (\dxg\z_{1}^{p}\dxg\z_{1}-\dxg\z_{2}^{p}\dxg\z_{2})\dxg W \right\}\|_{\infty} \nonumber \\
& & + C_{\ga} \| W\|_{C^{1}} \| \dxg\z_{1}^{p}\dxg\z_{1}-\dxg\z_{2}^{p}\dxg\z_{2}\|_{\infty} \label{eq:bitiyo}\\
& & + \left| \int_{\tor} W \left\{ \dxg\z_{1}^{p}\dxg\z_{1}-\dxg\z_{2}^{p}\dxg\z_{2}\right\}\right| \nonumber
\end{eqnarray} 
By the definition of $\dxi$ the first term on the right hand side of this equality is bounded by
\[ C_{\ga} \| W\|_{C^{1}} \| \dx\z_{1}^{p}\dx\z_{1}-\dx\z_{2}^{p}\dx\z_{2}\|_{\infty}
\]
which, by adding and subtracting the appropriate terms can be estimated by
\[ C_{\ga,\z_{1},\z_{2}} \| W\|_{C^{1}} \|\z_{1}-\z_{2}\|_{C^{1}}
\]
where $C_{\ga,\z_{1},\z_{2}}$ depends only on $C^{1}$ norms of $\ga,\gi,\z_{1}$ and $\z_{2}$.

Similarly we estimate the second summand in (\ref{eq:bitiyo}) by
\[ C_{\ga,\z_{1},\z_{2}} \| W\|_{C^{1}} \| \z_{1}-\z_{2}\|_{C^{1}}
\]
and the third summand by
\[ C_{\ga,\z_{1},\z_{2}} \| W\|_{\infty} \| \z_{1}-\z_{2}\|_{C^{1}}.
\]
Hence $\partial_{\ga}F_{(\ga,\z)}$ is continuous in $\z$.
The continuity of $\partial_{\ga}F_{(\ga,\z)}$ in $\z$ can be shown analogously. 
Therefore $(\z^{p},F(\ga,\z))$ defines a continuous differentiable map in a neighborhood of $(id,0)$.

\hfill $\Box$
\begin{rema}
There are different forms of the Hunter-Saxton equation considered in the literature. The equation introduced by J.K. Hunter and R. Saxton in \cite{HS} is
\begin{equation} 
\dx (\dt u +u\dx u)=\frac{1}{2}(\dx u)^{2}
\label{eq:HuS}
\end{equation}
whereas the bihamiltonian system studied in \cite{HuZh1} is the derivative of (\ref{eq:HuS}):
\begin{equation} 
\dx^{2} (\dt u +u\dx u)=\frac{1}{2}\dx((\dx u)^{2}).
\label{eq:HuZ}
\end{equation}
Note that $\frac{1}{2}\int_{\tor}|\dx u|^{2} dx=\frac{1}{2} \int_{\tor}|\dx u_{0}|^{2}dx$ is a conserved quantity. 
Then the local well posedness of the periodic Cauchy problem for (\ref{eq:HuZ}) with initial data $u(x,0)=u_{0}(x)$ is equivalent to that of 
\begin{eqnarray}
\dt u +u \dx u=\frac{1}{2}\dxi ((\dx u)^{2})+\bar{u}'
\label{eq:mean}
\end{eqnarray}
where $\bar{u}(t)=\int_{\tor}u \ dx $. 

On the other hand, for any given $\bar{u}(t)\in C^{1}([0,\infty), \mathbb{R})$ there is a unique solution $(\z, \ga)$ to the Cauchy problem for the equation
\[ \begin{array}{l}
\zd=\frac{1}{2} \dxig((\dxg (\z+\bar{u}))^{2}) \\
\gd=\z+\bar{u} \end{array}
\]
such that $u=\z\circ\gi +\bar{u}$ is a solution to the Cauchy problem for (\ref{eq:mean}) by a simple modification of theorem \ref{th:C1}. Therefore for any given $\bar{u}(t)\in C^{1}([0,\infty), \mathbb{R})$, there is a unique solution $u \in C^{0}([0,T), C^{1}(\tor))\cap C^{1} ([0, T), L^{\infty} (\tor))$ to the Cauchy problem for the equation (\ref{eq:HuZ}) such that $\int_{\tor}\dt u \ dx= \bar{u}'$ for initial data $u_{0}\in C^{1}(\tor)$.
\label{rem:1}
\end{rema}

\begin{rema}
In \cite{Yi} Z.Yin shows that the initial value problem for (\ref{eq:HuZ}) has a family of solutions in $\sob(\tor)$ for $s>3/2$. Our approach clarifies the uniqueness issue for this problem by identifying the relation between $\bar{u}'$ in (\ref{eq:mean}) and a solution to the initial value problem for (\ref{eq:HuZ}) (see remark \ref{rem:1}) as well as improving the well-posedness result from $\sob$ for $s>3/2$ to $C^{1}$. 
\label{rem:2}
\end{rema}

\begin{rema}
The modified Hunter-Saxton equation
\[ \dx (\dt u +u^{p}\dx u)= \frac{1}{2}(\dx u^{p}\dx u )
\]
(in the form it is introduced in \cite{HS}) can be treated similarly by reformulating the problem as the Cauchy problem for the equation
\[ \begin{array}{l}
\zd=\frac{1}{2} \dxig((\dxg (\z+\bar{u})^{p}\dxg \z ) \\
\gd=(\z+\bar{u})^{p}. \end{array}
\]
\end{rema}
\section{Well-posedness in $\sob$}

In proving theorem \ref{th:2} we implement the same approach of rewriting the Cauchy problem for (mHS) as an initial value problem for an ordinary differential equation on the product space $\dif \times\sob(\tor)$ where $\dif$ is the space of $\sob$ class diffeomorphisms on $\tor$.
Here we repeatedly use the Schauder ring property, Sobolev lemma and the following well-known property of Sobolev spaces that we refer to as the composition lemma.
\begin{lema}{\bf [Composition lemma]}
Let $s > 3/2 $ , $u \in \sob$ and $\ga$ be a $\sob$ class bijection from $\tor$ to $\tor$ whose inverse
$\gi$ is also of class $\sob$. Then, $u\circ \ga \in \sob$ and the estimate
\begin{equation}
\| u \circ \ga \|_{\sob} \leq C_{\z} (1 + \| \ga \|_{\sob}^{s}) \| u \|_{\sob}
\label{eq:com}
\end{equation}
holds where $C_{\ga}$ depends only on $C^{1}$ norms of $\ga$ and $\gi$.
\label{lem:com}
\end{lema}

The next proposition is the $\sob$ version of proposition \ref{pro:odeC1}. It is a restatement of the Cauchy problem for (mHS) with $\sob$ initial data as an initial value problem for an ordinary equation on $\dif\times\sob(\tor)$.
We first state the Cauchy problem for (mHS) on the product space $\dif\times\sob(\tor)$. 
\begin{prop}
A function $u \in \sod$ is a solution to the Cauchy problem
\[ (\dt u +u^{p} \dx u)=\frac{p}{2} \dxi (u^{p-1} (\dx u)^{2}) ; \ u(x,0)=u_{0}(x)
\]
if and only if  $u$ can be written as $\z \circ \gi$ where $(\z , \ga)\in\dif\times\sob$ is a solution to
\begin{equation}
\begin{array}{l} 
\gd = \z^{p}, \\ \\
\zd =F(\ga, \z) := \frac{\textstyle p}{\textstyle 2} ( \dxi ( (\z \circ \gi)^{p-1} (\dx (\z \circ \gi))^{2} ) ) \circ \ga.
\end{array}
\label{eq:**}
\end{equation}
with initial data $\z (x,0)=u_{0}(x)$ and $\ga (x,0) = id_{x} $.
\label{pro:ode}
\end{prop}

\pf Given a function $u \in \sod(\tor) $ which is a solution to the problem (\ref{eq:**}) there is a $\ga \in \dif $
satisfying
\[ \gd = u^{p} \circ \ga, \ \ga (x,0)= id_{x}
\]
(see \cite{BourBre}).
We set $\z = u \circ \ga $. Then we have 
\[ \zd = (\dt u + u^{p} \dx u) \circ \ga .\]
But $u$ is a solution to the Cauchy problem for (mHS) therefore (\ref{eq:**}) holds for $(\ga,\z)$. Conversely assume that we have a solution $(\ga,\z )\in \dif\times\sod(\tor) $ to (\ref{eq:**}) and let $u= \z\circ \gi$. Then \( \z=u \circ \ga \)
and therefore
\[ \dot{\z}=(\dt u+ (\gd\circ \gi) \dx u) \circ \ga .
\]
But since $(\ga,\z)$ is a solution to (\ref{eq:**}) $u$ is a solution to the Cauchy problem for
(mHS). 

\hfill $\Box$

The next lemma is used in the proof of theorem \ref{th:2} to show that $\ga\rightarrow\dxig$ is a continuous map from $\dif$ into $L(\sobb,\sob)$.
\begin{lema}
Let $s>3/2$. Then we have 
\[ \|\dxi(u\circ\gi)\circ\ga-\dxi u\|_{\sob}\leq C \|u\|_{\sobb} \| \ga - id_{x}\|_{\sob}(\| \ga \|_{\sob}+1)
\]
for $u\in\sobb(\tor)$ and $\ga \in\diff $.
\label{lema:lipsc} 
\end{lema} 
\pf Using the identity in (\ref{eq:identity}) we can write the $\sob$ norm to be estimated as follows:
\begin{equation*}
\|\dxi(u\circ\gi)\circ\ga-\dxi u\|_{\sob} = \| \dxi(u(\dx\ga -1))\|_{\sob}.
\end{equation*}
Therefore by Poincar\'{e} inequality, since $\dxi f$ has zero mean we have
\begin{equation}
\|\dxi(u\dx\ga-u)\|_{\sob}\leq \| u(\dx\ga -1) \|_{\sobb}.
\label{eq:pl}
\end{equation}
Finally using the Schauder ring property and Sobolev lemma with (\ref{eq:pl}) we obtain
\[ \|\dxi(u\circ\gi)\circ\ga-\dxi u\|_{\sob}\leq C \|u\|_{\sobb} \| \ga - id_{x}\|_{\sob}.
\]

\hfill $\Box$

Now we are ready to proceed to the proof of theorem \ref{th:2}. We will show that $F(\ga,\z)$ is a bounded map with bounded linear directional derivatives which are continuous. Then theorem \ref{th:2} follows from the fundamental theorem for ordinary differential equations on Banach spaces.

\noindent {\bf \textit{Proof of theorem \ref{th:2}. }} We start by showing that 
\[ F(\ga,\z)=\frac{1}{2}\dxig\{ \dxg\z^{p}\dxg\z\} \] 
defines a bounded map from $\dif\times\sob(\tor)$ into $\sob(\tor)$. By the composition lemma we have 
\begin{equation}
\| F(\ga,\z)\|_{\sob}\leq C_{\ga}\| \dxi (\dx (\z^{p}\circ\gi)\dx(\z\circ\gi))\|_{\sob}
\label{eq:yildiz}
\end{equation}
where $C_{\ga}$ depends only on $\sob$ norm of $\ga$ and $\gi$. Using Poincar\'{e} inequality we obtain the estimate
\[ \leq C^{\ga} \| \dx (\z^{p}\circ\gi)\dx(\z\circ\gi)\|_{\sobb}. 
\]
Then Schauder ring property with composition lemma gives the bound
\[ \leq C_{\ga} \| \z\|_{\sob}^{p+1}.
\]

Next we show that the directional derivatives
\begin{eqnarray}
\partial_{\ga}F_{(\ga,\z)}(W)=-\frac{1}{2}\dxig (\dxg \z^{p}\dxg \z \dxg W) \nonumber \\
+\frac{1}{2}\int_{\tor} W \dxg \z^{p}\dxg \z \ dx,
\label{eq:dgf1}  
\end{eqnarray}
\begin{equation}
\partial_{\z}F_{(\ga,\z)}(W)=\frac{1}{2}\dxig \{ p\dxg (\z^{p-1}W)\dxg \z + \dxg \z^{p}\dxg W \}
\label{eq:dzf1}  
\end{equation}
define bounded linear maps. We estimate the first summand on the right hand side of
\begin{eqnarray}
\|\partial_{\ga}F_{(\ga,\z)}(W)\|_{\sob} &\simeq & \|\dxig (\dxg \z^{p}\dxg \z \dxg W)\|_{\sob} \nonumber \\
& & +\left| \int_{\tor} W \dxg \z^{p}\dxg \z \ dx \right| 
\label{eq:ikiyildiz1}  
\end{eqnarray}
and $\| \partial_{\z}F_{(\ga,\z)}W\|_{\sob}$ using the Schauder ring property, Poincar\'{e} lemma and composition lemma by
\[ C_{\ga} \| W\|_{\sob}\| \z\|_{\sob}^{p} (1+\|\z\|_{\sob})
\]
as we estimated (\ref{eq:yildiz}) above.

For the absolute value term in (\ref{eq:ikiyildiz1}) we apply Sobolev lemma to obtain the bound
\[ C_{\ga} \| W\|_{\sobb}\|\z\|_{\sob}^{p+1}.
\]
Hence the directional derivatives $\partial_{\ga}F_{(\ga,\z)}$ and $\partial_{\z}F_{(\ga,\z)}$ define bounded linear maps.
The continuity of the directional derivatives can be proved as follows.

We observe that (like in the previous section) to prove the continuity of $\ga\rightarrow\partial F_{(\ga,\z)}$ it suffices to prove the continuity of the maps
\begin{eqnarray}
 & & \ga\rightarrow\partial_{\ga}\dxig\in L(\sob,L(\sobb,\sob)) \label{eq:onuc} \\
 & & \ga\rightarrow\dxig\in L(\sobb,\sob) \label{eq:ondort} \\
 & & \ga\rightarrow\dxg\in L(\sob,\sobb) \label{eq:onbes}
\end{eqnarray}
in a neighborhood of $(id,0)\in \dif\times \sob(\tor)$.

The identity in (\ref{eq:identity}), Poincar\'{e} inequality and Schauder ring property lead to the following estimate for the map in (\ref{eq:ondort})
\[ \| \dxig f - \dxi f\|_{\sob} \leq C_{\ga} \| f\|_{\sobb} \| \ga -id \|_{\sob}
\]
and hence to the continuity of this map.

In order to show that the map in (\ref{eq:onuc}) is continuous it is sufficient to obtain an appropriate estimate on the norm 
\begin{eqnarray*}
\lefteqn{\| \partial_{\ga}\dxig W- \partial_{\ga} (\dxi)_{id}W\|_{\sob}} \\
&\leq & \| \dxig \{ \dxg\z^{p}\dxg\z\dxg W\}- \dxi (\dx\z^{p}\dx\z\dx W)\|_{\sob} \\
& & + \left| \int_{\tor} W \{ \dxg \z^{p}\dxg \z - \dx \z^{p} \dx \z \} dx\right|. 
\end{eqnarray*}
The second term on the right hand side above is bounded using Sobolev lemma and the Schauder ring property by
\[ C_{\ga} \| W\|_{\sobb} \| \z\|_{\sob}^{p+1} \|\ga -id\|_{\sob}.
\]
Adding and subtracting the appropriate terms and using Poincar\'{e} inequality we obtain the following estimate for the first term 
\begin{eqnarray}
\lefteqn{\| \dxig \{ \dxg\z^{p}\dxg\z\dxg W\}- \dxi (\dx\z^{p}\dx\z\dx W)\|_{\sob}}  \nonumber \\
&\leq & \| \dxig (\dxg\z^{p}\dxg\z\dxg W)- \dxig (\dx\z^{p}\dx\z\dx W)\|_{\sob} \label{eq:deli}\\
& & + \| \dxig (\dx\z^{p}\dx\z\dx W)-\dxi (\dx\z^{p}\dx\z\dx W)\|_{\sob} \label{eq:kupeli}.
\end{eqnarray}
For (\ref{eq:kupeli}) lemma \ref{lema:lipsc} and the Schauder ring property give us the following estimate
\[ C_{\ga} \|\z\|_{\sob}^{p+1} \|W\|_{\sob} \|\ga -id\|_{\sob}.
\]
Applying composition lemma and Poincar\'{e} inequality to (\ref{eq:deli}) we obtain
\[
\| \dxig (\dxg\z^{p}\dxg\z\dxg W)- \dxig (\dx\z^{p}\dx\z\dx W)\|_{\sob} \]
\[ \leq  C_{\ga} \|\z\|_{\sob}^{p+1} \|W\|_{\sob} \|\ga -id\|_{\sob}.\]
Therefore the map in (\ref{eq:onuc}) is continuous.

It is not difficult to obtain the following inequality which implies that the map in (\ref{eq:onbes}) is continuous using  Schauder ring property    
\[ \| \dxg W - \dx W\|_{\sobb} \leq C_{\ga} \| W\|_{\sob} \|\ga -id\|_{\sob}.
\]

In order to conclude that $F(\ga,\z)$ defines a continuously differentiable map $\dif\times\sob$ into $\sob$ it is sufficient to observe that both directional derivatives  $\partial_{\ga}F_{(\ga,\z)}$ and $\partial_{\z}F_{(\ga,\z)}$ are continuous in $\z$ as well. This completes the proof of theorem \ref{th:2}.

\hfill $\Box$

\section{Analytic regularity}

In this section we give a proof of theorem \ref{th:analy_reg} that states the analytic regularity (i.e., 
existence and uniqueness of analytic solutions for analytic initial data) of the Cauchy problem for (mHS).

The classical Cauchy-Kowalevski theorem does not apply to the modified Hunter-Saxton equation (mHS). However a contraction argument on a scale of Banach spaces can be used for the nonlocal form (mHS) of this equation. The conditions under which such a contraction argument can be applied to prove an analytic regularity result are given in \cite{Nis} and  \cite{BG} in the form of an abstract Cauchy-Kowalevski theorem.    
The following scale of Banach spaces is appropriate for our problem.
\begin{defi}
Let $s$ be a positive real number. The collection $\{ X_{s}\}_{s>0}$ of Banach spaces $X_{s}$ is called a decreasing scale
of Banach spaces if $s'<s$ implies $X_{s}\subset X_{s'}$ and $|||\cdot |||_{s'}\leq |||\cdot |||_{s}$.
\end{defi} 
For $s>0$, let the spaces $E_{s}$ be defined as
\[ E_{s}=\left\{ u \in C^{\infty}(\tor) : \int_{\tor} u \ dx=0 \ and \ ||| u |||_{s}=\sup_{k\geq 0} \frac{\| \dx^{k}.
u \|_{\sob}s^{k}}{k! / (k+1)^{2}} < \infty \right\},
\]

We use decreasing scales of Banach spaces $X_{s}:= E_{s}\times E_{s}$ given by Cartesian products
of the spaces $E_{s}$ for the contraction argument. The norm $|||\cdot |||_{X_{s}}$ can be chosen to be any of the
standard product norms on $E_{s}\times E_{s}$. An important result for the spaces $E_{s}$ is the following algebra property.   
\begin{lema}
Let $0<s<1$. There is a constant $c>0$ which is independent of $s$ such that we have
\[ ||| uv |||_{s}\leq c |||u|||_{s} |||v|||_{s}
\]
for any $u,v \in E_{s}$.
\label{lema:product}
\end{lema}

We do not repeat the proof of this lemma here but refer to \cite{HM1} where A.Himonas and G.Misio\l ek prove the analytic regularity of solutions to the Cauchy problem for the Camassa-Holm equation.

First we rewrite the equation (mHS) in a more convenient form.
Let the operators $P_{1}$ and $P_{2}$ be defined as
\begin{align*} 
& P_{1}(u):=-\dx u, \\ 
& P_{2}(u):=\dxi u
\end{align*} 
and let $u_{1}:=u, u_{2}:=\dx u$. Then we can write the equation
(mHS) in terms of $(u_{1},u_{2})$ as follows 
\begin{equation} \begin{array}{l}
\dt u_{1}=F_{1}(u_{1},u_{2}):=\frac{\textstyle 1}{\textstyle p+1} P_{1}(u_{1}^{p+1})+\frac{\textstyle p}{\textstyle
2}P_{2}(u_{1}^{p-1}u_{2}^{2})
\\ 
\dt u_{2}=F_{2}(u_{1},u_{2}):=P_{1}(u_{1}^{p}u_{2})+\frac{\textstyle p}{\textstyle 2}u_{1}^{p-1}u_{2}^{2}
\end{array}
\label{eq:syst}
\end{equation}

The following two lemmas give suitable bounds on the operators $P_{1}$ and $P_{2}$ to prove theorem
\ref{th:analy_reg}. 

\begin{lema}
For $0<s'<s<1$, the estimate
\[ ||| P_{1}(u)|||_{s'} \leq \frac{1}{s-s'} ||| u |||_{s} \]
holds.
\label{lema:P1}
\end{lema}

\begin{lema}
For any $0<s<1$, the estimate 
\[ ||| P_{2}(u)|||_{s} \leq ||| u |||_{s} \]
holds.
\label{lema:P2}
\end{lema}

Now we are ready to prove Theorem \ref{th:analy_reg}. We will show that all three conditions of the abstract version of the Cauchy-Kowalevski theorem that we include in the appendix for the convenience of the reader as it is stated in \cite{Nis}, hold for (mHS) on the scale $\{ X_{s}\}_{0<s<1}$.

{\bf Proof of Theorem \ref{th:analy_reg}.} Let $u=(u_{1},u_{2})$ and $F=(F_{1},F_{2})$ in (\ref{eq:syst})
and let $X_{s}$ be a decreasing scale of Banach spaces defined as $X_{s}=E_{s}\times E_{s}$. Then we only need to verify
the first two conditions of the abstract Cauchy-Kowalevski theorem since the map $F(u_{1},u_{2})$ does not depend on $t$
explicitly.

Clearly, $t\longmapsto F(t,u(t))=(F_{1}(u_{1},u_{2}),F_{2}(u_{1},u_{2}))$ is holomorphic if $t\longmapsto u_{1}(t)$ and 
$t\longmapsto u_{2}(t)$ are both holomorphic. Therefore, to verify the first condition of the abstract theorem, we only
need to show that for $s'<s$, 
$F_{1}(u_{1},u_{2})$ and
$F_{2}(u_{1},u_{2})$  are in $ E_{s'}$ if $u_{1},u_{2}\in E_{s}$. We begin with estimates on $F_{1}$:
\[ ||| F_{1}(u_{1},u_{2}) |||_{s'}= ||| \frac{1}{p+1} P_{1}(u_{1}^{p+1})
+\frac{p}{2}P_{2}(u_{1}^{p-1}(u_{2})^{2}) |||_{s'}
\]
\[ \leq c_{p}(||| P_{1}(u_{1}^{p+1}) |||_{s'}+||| P_{2}(u_{1}^{p-1}(u_{2})^{2})  |||_{s'}).
\]
By Lemma \ref{lema:P1} and Lemma \ref{lema:P2} we have the following bound on this last term 
\[ \leq c_{p}\left( \frac{1}{s-s'}||| u_{1}^{p+1} |||_{s}+ |||u_{1}^{p-1}(u_{2})^{2} |||_{s'}\right).
\]
Then by Lemma \ref{lema:product} we obtain the estimate
\[ \leq \frac{c_{p}}{s-s'}||| u_{1} |||_{s}^{p+1}+ c_{p}|||u_{1}|||_{s'}^{p-1}|||u_{2}
|||_{s'}^{2}.
\]
Similarly, for $F_{2}$, we have
\begin{align*}
||| F_{2}(u_{1},u_{2}) |||_{s'} & =||| P_{1}(u_{1}^{p}u_{2})+\frac{p}{2}u_{1}^{p-1}u_{2}^{2} |||_{s'} \\
& \leq |||  P_{1}(u_{1}^{p}u_{2}) |||_{s'}+\frac{p}{2}|||u_{1}^{p-1}u_{2}^{2} |||_{s'} \\
& \leq c \left( \frac{1}{s-s'}|||u_{1}^{p}u_{2} |||_{s}+\frac{p}{2}|||u_{1}|||_{s'}^{p-1}|||u_{2}|||_{s'}^{2}\right) \\
& \leq \frac{c}{s-s'}|||u_{1}|||_{s'}^{p}|||u_{2}|||_{s}+c \frac{p}{2}|||u_{1}|||_{s'}^{p-1}|||u_{2}|||_{s'}^{2}.
\end{align*}
We proceed to establish the second condition of the abstract Cauchy-Kowalevski theorem. First we show that for some $c$
independent of $t$ the estimates
\begin{equation}  
||| F_{1}(u_{1},u_{2})- F_{1}(v_{1},v_{2})|||_{s'} \leq \frac{c}{s-s'} |||u-v|||_{X_{s}}
\label{est1}
\end{equation}
and 
\begin{equation}  
||| F_{2}(u_{1},u_{2})- F_{2}(v_{1},v_{2})|||_{s'} \leq \frac{c}{s-s'} |||u-v|||_{X_{s}}
\label{est2}
\end{equation}
hold.

In order o obtain the estimate (\ref{est1}) we use the triangle inequality and Lemma \ref{lema:P1} with Lemma \ref{lema:P2}:
\begin{align*}
||| F_{1}(u_{1},u_{2})- F_{1}(v_{1},v_{2})|||_{s'}& \leq c_{p} |||P_{1}(u_{1}^{p+1}) + P_{2}(u_{1}^{p-1}u_{2}^{2}) -P_{1}(v_{1}^{p+1}) -
P_{2}(v_{1}^{p-1}v_{2}^{2}))|||_{s'} \\
& \leq c_{p} (||| P_{1}(u_{1}^{p+1}-v_{1}^{p+1}) |||_{s'} + ||| P_{2}(u_{1}^{p-1}u_{2}^{2}-v_{1}^{p-1}v_{2}^{2})|||_{s'}) \\
& \leq \frac{c_{p}}{s-s'} (||| (u_{1}^{p+1}-v_{1}^{p+1}) |||_{s} + ||| (u_{1}^{p-1}u_{2}^{2}-
v_{1}^{p-1}v_{2}^{2})|||_{s'})  
\end{align*}
Then we add and subtract $u_{1}^{p-1}v_{2}^{2}$ and use Lemma \ref{lema:product} to get
\begin{align*}
& = \frac{c_{p}}{s-s'} (||| (u_{1}^{p+1}-v_{1}^{p+1}) |||_{s} + ||| u_{1}^{p-1}(u_{2}^{2}-v_{2}^{2})+
(v_{1}^{p-1}-u_{1}^{p-1})v_{2}^{2}|||_{s'}) \\
& \leq \frac{c_{p}}{s-s'} (||| u_{1}^{p+1}-v_{1}^{p+1} |||_{s} + ||| u_{1}^{p-1} |||_{s'} ||| u_{2}^{2}-v_{2}^{2}
|||_{s'} + |||v_{1}^{p-1}-u_{1}^{p-1} |||_{s'} ||| v_{2}^{2}|||_{s'}).
\end{align*}
Note that using the identity
\[ u_{1}^{p}-v_{1}^{p}=(u-v)(u^{p-1}+u^{p-2}v+u^{p-3}v^{2}+...+uv^{p-2}+v^{p-1}).
\]
and the assumption that $|||u|||_{s}<R$ and $ |||v|||_{s}<R$ we have
\begin{align*}
||| F_{1}(u_{1},u_{2})- F_{1}(v_{1},v_{2})|||_{s'}  \leq & \frac{c_{p}}{s-s'} \left((p+1)R^{p} 
||| u_{1}-v_{1} |||_{s} \right) \\
& + 2R^{p} ||| u_{2}-v_{2} |||_{s'} + R^{2}(p-1)R^{p-2} |||v_{1}-u_{1} |||_{s'} \\ 
\leq & \frac{c_{p,R}}{s-s'} |||u-v|||_{X_{s}}
\end{align*}
where $u=(u_{1},u_{2})$ and $v=(v_{1},v_{2})$.

The estimate (\ref{est2}) we first use triangle inequality
\begin{align*}
|||F_{2}(u_{1},u_{2})-F_{2}(v_{1},v_{2}) |||_{s'} & =
||| P_{1}(u_{1}^{p}u_{2}) + \frac{p}{2}u_{1}^{p-1}u_{2}^{2}- P_{1}(v_{1}^{p}v_{2}) -\frac{p}{2}v_{1}^{p-1}v_{2}^{2}|||_{s'}
\\
& \leq c_{p} (|||P_{1}(u_{1}^{p}u_{2}-v_{1}^{p}v_{2}) |||_{s'} + |||u_{1}^{p-1}u_{2}^{2}-v_{1}^{p-1}v_{2}^{2}|||_{s'}).
\end{align*}
Then using Lemma \ref{lema:P1} we estimate this last term by
\[ = \frac{c_{p}}{s-s'} (||| u_{1}^{p}u_{2}-v_{1}^{p}v_{2} |||_{s} +
|||u_{1}^{p-1}u_{2}^{2}-v_{1}^{p-1}v_{2}^{2}|||_{s'}).
\]
Adding and subtracting $ u_{1}^{p}v_{2}$ and $ u_{1}^{p-1}v_{2}^{2} $ and using Lemma \ref{lema:product} we obtain
\begin{align*}
= \frac{c_{p}}{s-s'} & (||| u_{1}^{p}(u_{2}-v_{2})+ (u_{1}^{p}-v_{1}^{p})v_{2} |||_{s} +
|||u_{1}^{p-1}(u_{2}^{2}-v_{2}^{2}) + (u_{1}^{p-1}-v_{1}^{p-1})v_{2}^{2}|||_{s'}) \\
\leq \frac{c_{p}}{s-s'} & (||| u_{1}^{p}(u_{2}-v_{2}) |||_{s}+ ||| (u_{1}^{p}-v_{1}^{p})v_{2} |||_{s} +
|||u_{1}^{p-1}(u_{2}^{2}-v_{2}^{2}) |||_{s'} \\
& + ||| (u_{1}^{p-1}-v_{1}^{p-1})v_{2}^{2}|||_{s'}) \\
\leq \frac{c_{p}}{s-s'} & (||| u_{1}^{p} |||_{s} |||u_{2}-v_{2} |||_{s}+ ||| u_{1}^{p}-v_{1}^{p} |||_{s} ||| v_{2}
|||_{s} + |||u_{1}^{p-1} |||_{s'} |||u_{2}^{2}-v_{2}^{2} |||_{s'}  \\
& + ||| u_{1}^{p-1}-v_{1}^{p-1}|||_{s'} |||v_{2}^{2}|||_{s'}).
\end{align*}
Again assuming that $|||u|||_{s}<R$ and $ |||v|||_{s}<R$ we obtain 
\[ |||F_{2}(u_{1},u_{2})-F_{2}(v_{1},v_{2}) |||_{s'} \leq \frac{c_{p,R}}{s-s'}|||u-v|||_{X_{s}}.
\]
which completes the proof of Theorem \ref{th:analy_reg}.

\hfill $\Box$

\section*{Appendix An abstract Cauchy-Kowalevski theorem}

The classical Cauchy-Kowalevski theorem is an important tool in the study of
partial differential equations. An abstract version of this theorem which can be
applied to differential equations involving pseudodifferential operators has been
developed since '60s.

The general framework in the case of linear equations has appeared in the work of T. Yamanaka
\cite{Y} and L.V. Ovsjannikov \cite{Ovs1}. The first nonlinear form of the abstract
theorem presented by F. Treves \cite{Tre} has been improved by L. Nirenberg
\cite{Nir} and Nirenberg's version has been further simplified by T. Nishida
\cite{Nis}. We will use the following version of the abstract Cauchy-Kowalevski
theorem from \cite{Nis}. 

\begin{teor}
Consider the Cauchy problem 
\begin{equation}
\begin{array}{l}
\dt u = F(t, u(t)) \\ \\
 u(0)=0.
\end{array}
\label{eq:analytic}
\end{equation}
Let $\{ X_{s}\}_{0<s<1}$ be a scale of decreasing Banach spaces, so that for any $s'<s$ we have $X_{s}\subset
X_{s'}$ and $ ||| . |||_{s'} \leq ||| . |||_{s} $. 
Let $T,R$ and $C$ be positive numbers and suppose that $F$ satisfies the following conditions

1.) If for $0<s'<s<1$ the function $t\longmapsto u(t)$ is holomorphic in $|t|<T$ and continuous on $|t|\leq
T$ with values in $X_{s}$ and 
\[ \sup_{|t|\leq T} ||| u(t) |||_{s} <R,
\]
then $t\longmapsto F(t,u(t))$ is a holomorphic function on $|t|<T$ with values in $X_{s'}$.

2.) For any $0<s'<s\leq 1$ and any $u,v \in X_{s}$ with $|||u|||_{s}<R, ||| v |||_{s}<R$,
\[ \sup_{|t|\leq T} ||| F(t,u)-F(t,v) |||_{s'} \leq \frac{C}{s-s'}|||u-v |||_{s}.
\]

3.) There exists $M>0$ such that for, any $0<s<1$, 
\[ \sup_{|t|\leq T} |||F(t,0) |||_{s} \leq \frac{M}{1-s} .
\]
Then there exists a $T_{0}\in (0,T)$ and a unique function $u(t)$, which for every $s \in (0,1)$ is
holomorphic in $|t|<(1-s)T_{0}$ with values in $X_{s}$, and is a solution to the initial value problem 
(\ref{eq:analytic}).
\label{th:cau_kowa}
\end{teor}

\vspace{.5in}

{\em Feride Tiglay

University of New Orleans,
 
Department of Mathematics,

Lakefront,

New Orleans LA 70148 }


\begin{thebibliography}{HuSa}

\bibitem[AK]{ArKh} V.I. Arnold and B.A. Khesin, {\em Topological Methods in Hydrodynamics}, Springer
Verlag, New York 1998.

\bibitem[Arn]{Arn} V. Arnold, {\em Sur la g\'{e}om\'{e}trie diff\'{e}rentielle des
groupes de Lie de dimension infinie et ses applications \`{a} l'hydrodynamique des fluides
parfaits}, Ann. Inst. Fourier, Grenoble {\bf 16} (1966).

\bibitem[BB]{BourBre} J.P. Bourguignon and H. Brezis, {\em Remarks on the Euler equation}, J.
Functional Analysis {\bf 15} (1974).

\bibitem[BG]{BG}M.S. Baouendi and C. Goulaouic, {\em Sharp estimates for analytic pseudodifferential
operators and application to Cauchy problems}, J. Diff. Eq. {\bf 48} (1983).

\bibitem[BSS1]{BSS1} R. Beals, D. Sattinger and J. Szmigielski,{\em Acoustic scattering and the extended Korteweg-de Vries equation,}, Advances in Mathematics, {\bf 140} (1998).

\bibitem[BSS2]{BSS2} R. Beals, D.H. Sattinger and J. Szmigielski, {\em Inverse scattering solutions of the Hunter-Saxton equation}, Appl. Anal. {\bf 78}  (2001).

\bibitem[Di]{Dieud} J. Dieudonn\'{e}, {\em Foundations of Modern Analysis}, Academic Press, New York
and London, 1960.

\bibitem[Eb1]{Eb1} D.G. Ebin, {\em Espace des m\'{e}triques riemanniennes et mouvement des fluides
via les vari\'{e}t\'{e}s d'applications}, Centre de Math\'{e}matiques de l'Ecole Polytechnique et
Universit\'{e} Paris {\rm 8} 1972.

\bibitem[Ee]{Eells} J. Eells, {\em A setting for global analysis}, Bull. Amer. Math. Soc. {\bf 72}
(1966).

\bibitem[EMa]{EMa} D.G. Ebin and J. Marsden, {\em Groups of diffeomorphisms and the motion of an
incompressible fluid}, Ann. Math.
{\bf 92} (1970).

\bibitem[HM1]{HM1}A.A. Himonas and G. Misio\l ek, {\em Remarks on an integrable evolution equation},
Geometry
and analysis on finite and infinite dimensional Lie groups, Banach Center Publ. {\bf 55}, Polish
Acad. Sci.,
Warsaw (2002).

\bibitem[HM2]{HM2} A.A. Himonas and G. Misio\l ek, {\em The Cauchy problem for an integrable
shallow-water equation}, Differential and Integral Equations {\bf 14} (2001).

\bibitem[HuSa]{HS} J.K. Hunter and R. Saxton, {\em Dynamics of director fields}, SIAM J.
Appl. Math. {\bf 51}, No. 6 (1991). 

\bibitem[HZ1]{HuZh1} J.K. Hunter and Y. Zheng, {\em On a completely integrable nonlinear hyperbolic
variational equation}, Physica D {\bf 79}, 1994.

\bibitem[K]{Kato} T. Kato, {\em Liapunov functions and monotonicity in the Navier-Stokes equation},
Lecture Notes in Math., 1450, 
Springer, Berlin 1990.

\bibitem[KaM]{KaM}T. Kato and K. Masuda, {\em Nonlinear evolution equations and analyticity I}, Ann. de l'Inst. H. Poincar\'{e} {\bf 3} (1986).

\bibitem[KM]{KhMi} B.A. Khesin and G. Misio\l ek, {\em Euler equations on homogeneous spaces and
Virasoro orbits}, Adv. Math.   {\bf 176}, 2003.

\bibitem[Mis1]{Mis1} G. Misio\l ek, {\em Classical solutions of the periodic Camassa-Holm equation},
GAFA {\bf 12} (2002).

\bibitem[Mis3]{Mis3} G. Misio\l ek, {\em A shallow water equation as a geodesic flow on the
Bott-Virasoro group}, Journal of Geometry and Physics {\bf 24} (1996).

\bibitem[Nir]{Nir}L. Nirenberg, {\em An abstract form of the nonlinear Cauchy-Kowalevski theorem},
J. Differential Geometry {\bf 6} (1972).

\bibitem[Nis]{Nis} T. Nishida, {\em A note on a theorem of Nirenberg}, J. Differential Geometry {\bf
12} (1977).

\bibitem[OK]{OK} V. Ovsienko and B. Khesin, {\em Korteweg-de Vries superequations as an Euler
equation}, Functional Anal. Appl. {\bf 21} (1987).

\bibitem[Ovs1]{Ovs1}L.V. Ovsjannikov, {\em A singular operator in a scale of Banach spaces},
Dokl. Akad. Nauk SSSR, {\bf 163} (1965).

\bibitem[Ovs2]{Ovs2}L.V. Ovsjannikov, {\em A nonlinear Cauchy problem in a scale of Banach spaces},
Dokl. Akad. Nauk SSSR, {\bf 200} (1971); Soviet Math. Dokl. {\bf 12} (1971).

\bibitem[R]{R}E.G. Reyes, {\em The soliton content of the Camassa-Holm and Hunter-Saxton equations}, Proc. Inst. Math. NAS**** {\bf 43} (2002).


\bibitem[Sat]{Sattinger} D.H. Sattinger, {\em Scaling, mathematical modeling, \& Integrable systems}, DMV-Seminar Series, Band 28,  Birkhauser 1998.

\bibitem[Tay1]{Tay1} M.E. Taylor, {\em Pseudodifferential Operators and Nonlinear PDE},
Birkh\"{a}user Boston 1991.

\bibitem[Tay2]{Tay2} M.E. Taylor, {\em Finite and infinite dimensional Lie groups and evolution
equations}, Classroom Notes, University of North Carolina-Chapel Hill, Spring 2003.

\bibitem[Tre]{Tre}F. Treves, {\em An abstract nonlinear Cauchy-Kovalevska theorem}, Trans. Amer. 
Math. Soc. {\bf 150}
(1970).

\bibitem[Tru]{Tru}E. Trubovitz, {\em The inverse problem for periodic potentials}, Comm. Pure Appl. Math.  {\bf 30}
(1977).

\bibitem[Ya]{Y}T. Yamanaka, {\em Note on Kowalevskaja's system of partial differential equations},
Comment. Math. Univ. St. Paul. {\bf 9} (1960).

\bibitem[Yi]{Yi}Z. Yin, {\em On the structure of solutions to the periodic Hunter-Saxton equation},
SIAM J. Math. Anal. {\bf 36} (2004).

\end{thebibliography}
\end{document}